\documentclass[11pt]{amsart}
\usepackage{amsmath,amsfonts,amsthm,amsopn,amssymb}
\usepackage{cite,marginnote}
\usepackage{amsmath}
\usepackage{color,enumitem,graphicx}
\usepackage[colorlinks=true,urlcolor=blue,
citecolor=red,linkcolor=blue,linktocpage,pdfpagelabels,
bookmarksnumbered,bookmarksopen]{hyperref}
\usepackage[english]{babel}
\usepackage[left=2.8cm,right=2.8cm,top=3cm,bottom=3cm]{geometry}

\usepackage[hyperpageref]{backref}

\numberwithin{equation}{section}

\makeindex
\newtheorem{theorem}{Theorem}[section]
\newtheorem{proposition}[theorem]{Proposition}
\newtheorem{lemma}[theorem]{Lemma}
\newtheorem{remark}[theorem]{Remark}
\newtheorem{example}[theorem]{Example}
\newtheorem{corollary}[theorem]{Corollary}
\newtheorem{definition}[theorem]{Definition}

\newcommand{\RN}{\mathbb R^N}

\newcommand{\iy}{\infty}

\newcommand{\s}{\section}

\newcommand{\DD}{\Delta}
\newcommand{\g}{\gamma}

\newcommand{\R}{\mathbb R}
\newcommand{\al}{\alpha}

\newcommand{\re}[1]{(\ref{#1})}

\newcommand{\rg}{\rightarrow}

\newcommand{\e}{\varepsilon}

\newcommand{\bt}{\begin{theorem}}
\newcommand{\et}{\end{theorem}}
\newcommand{\bl}{\begin{lemma}}
\newcommand{\el}{\end{lemma}}
\newcommand{\bd}{\begin{definition}}
\newcommand{\ed}{\end{definition}}
\newcommand{\bc}{\begin{corollary}}
\newcommand{\ec}{\end{corollary}}
\newcommand{\bp}{\begin{proof}}
\newcommand{\ep}{\end{proof}}
\newcommand{\bx}{\begin{example}}
\newcommand{\ex}{\end{example}}
\newcommand{\bi}{\begin{exercise}}
\newcommand{\ei}{\end{exercise}}
\newcommand{\bo}{\begin{proposition}}
\newcommand{\eo}{\end{proposition}}
\newcommand{\br}{\begin{remark}}
\newcommand{\er}{\end{remark}}
\newcommand{\be}{\begin{equation}}
\newcommand{\ee}{\end{equation}}
\newcommand{\ba}{\begin{align}}
\newcommand{\ea}{\end{align}}
\newcommand{\bn}{\begin{enumerate}}
\newcommand{\en}{\end{enumerate}}
\newcommand{\bg}{\begin{align*}}
\newcommand{\eg}{\end{align*}}
\newcommand{\bcs}{\begin{cases}}
\newcommand{\ecs}{\end{cases}}

\newcommand{\bean}{\begin{eqnarray*}}
\newcommand{\eean}{\end{eqnarray*}}

\title[Singularly perturbed fractional Schr\"{o}dinger equation]{Singularly perturbed fractional Schr\"{o}dinger equation involving a general critical nonlinearity}

\author[H.\ Jin]{Hua Jin}
\author[W. B.\ Liu]{Wenbin Liu}
\author[J. J.\ Zhang]{Jianjun Zhang}

\address[H.\ Jin]{\newline\indent College of Science
\newline\indent
China University of Mining and Technology
\newline\indent
Xuzhou, 221116, China}
\email{\href{mailto:huajin@cumt.edu.cn}{huajin@cumt.edu.cn}}

\address[W. B.\ Liu]{\newline\indent College of Science
\newline\indent
China University of Mining and Technology
\newline\indent
Xuzhou, 221116, China}
\email{\href{mailto:liuwenbin-xz@163.com}{liuwenbin-xz@163.com}}

\address[J. J.\ Zhang]{\newline\indent College of Mathematica and Statistics
\newline\indent
Chongqing Jiaotong University
\newline\indent
Chongqing, 400074, China}
\email{\href{mailto:zhangjianjun09@tsinghua.org.cn}{zhangjianjun09@tsinghua.org.cn}}

\thanks{W. Liu is the corresponding author. }

\date{\today}

\begin{document}

\begin{abstract}
In this paper, we are concerned with the existence and concentration phenomena of solutions for the following singularly perturbed fractional Schr\"{o}dinger problem
\begin{align*}
\varepsilon^{2s}(-\Delta)^su+V(x)u=f(u) \ \ \ \mbox{in} \ \ \ \mathbb{R}^N,
\end{align*}
where $N>2s$ and the nonlinearity $f$ has critical growth. By using the variational approach, we construct a localized bound-state solution concentrating around an isolated component of the positive minimum point of $V$ as $\varepsilon\rightarrow 0$. Our result improves the study made in X. He and W. Zou ({\it Calc. Var. Partial Differential Equations}. 55-91(2016)), in the sense that, in the present paper, the {\it Ambrosetti-Rabinowitz} condition and {\it monotonicity} condition on $f(t)/t$ are not required.
\end{abstract}

\maketitle

\s{Introduction}
\renewcommand{\theequation}{1.\arabic{equation}}
\noindent {\bf 1.1 Background.}
In this paper, we are concerned with the standing waves for the nonlinear  fractional Schr\"{o}dinger equation
\begin{align}
\label{question1}
\varepsilon^{2s}(-\Delta)^su+V(x)u=f(u) \ \ \ \mbox{in} \ \ \mathbb{R}^N,
\end{align}
which is derived from the nonlinear fractional Schr\"{o}dinger equation
\begin{align}
\label{question2}
i\hbar \varphi_t-\hbar^2(-\Delta)^s\varphi-V(x)\varphi+f(\varphi)=0,\ \ \ (x,t)\in \mathbb{R}^N\times\mathbb{R},
\end{align}
where $\hbar$ is the Plank constant which is a very small physical quantity, $i$ is the imaginary unit and $N>2s$. The solutions of (\ref{question2}) with the form
\begin{align}
\label{standing}
\varphi(x,t)=e^{-iwt/\hbar}u(x), \ \ w\in \mathbb{R}.
\end{align}
is called standing waves. Assuming that $f:\mathbb{C}\rightarrow \mathbb{C}$ is continuous such that $f(e^{i\theta}u)=e^{i\theta}f(u)(u,\theta\in \mathbb{R}),$
and inserting (\ref{standing}) to (\ref{question2}), we have
\begin{align*}
\hbar^2(-\Delta)^su+(V(x)-w)u=f(u), \ \ \ x\in \mathbb{R}^N.
\end{align*}
Let $\hbar=\varepsilon^s$ and write $V-w$ as $V$, then we get (\ref{question1}). In the quantum mechanics, these standing waves are referred
to as semi-classical states, whose existence and concentration phenomena are particularly important as $\varepsilon\rightarrow 0$ .
$(-\Delta)^s(0<s<1)$ is the fractional Laplacian operator which can be seen as the infinitesimal generators of L\'{e}vy stable diffusion processes (see\cite {Applebaum}). This operator arises in many areas such as physics, biology, chemistry and finance (see\cite {Applebaum, laskin}). The fractional Schr\"{o}dinger equation is a fundamental equation of the fractional quantum mechanics. It was discovered by Laskin \cite{laskin1,laskin} as a result of extending the Feynman path integral, from the Brownian-like quantum mechanical path to L\'{e}vy-like one, where the Feynman path integral leads to the classical Schr\"{o}dinger equation, and the path integral over L\'{e}vy trajectories leads to the fractional Schr\"{o}dinger equation. For the further background in this field, we refer to \cite{EGE} and the references therein.

\

\noindent {\bf 1.2 Motivation.} An interesting class of solutions to (\ref{question1}) is a family of solutions that develop a spike shape around some point in $\mathbb{R}^N$ as $\varepsilon\rightarrow 0$. When $s=1$, equation (\ref{question1}) is reduced to a local elliptic equation
\begin{align}
\label{localquestion}
-\varepsilon^{2}\Delta u+V(x)u=f(u), \ \  x\in\mathbb{R}^N.
\end{align}
During the last three decades, many papers have been devoted to the singularly perturbed Schr\"{o}dinger equation (\ref{localquestion}) involving  subcritical growth or critical growth. Based on a Lyapunov-Schmidt reduction method, A. Floer and A. Weinstein \cite{A.Floer} first studied the existence of single-peak solutions for $N=1$ and $f(t)=t^3$. They constructed a single-peak solution which concentrates around any given nondegenerate critical point of $V$. More related results can be seen in \cite{J.Byeon1,M.delPino,weimingni,changfenggui,P.H.Rabinowitz} and the references therein. In these papers, the conditions such as Ambrosetti-Rabinowitz condition, monotonicity condition or nondegenerate condition are needed. To remove or weaken these conditions, J. Byeon and L. Jeanjean \cite{Byeon1} introduced a new penalization approach. Under the Berestycki-Lions conditions (see\cite{Berestycki}), the authors proved that for small $\varepsilon>0$, there exists a positive solution which clusters near a local minimum point of $V$ if $V$ satisfies

$(V_1)$ $V\in C(\mathbb{R}^N,\mathbb{R})$ and $0<V_0=\inf_{x\in \mathbb{R}^N}V(x)$,

$(V_2)$ there is a bounded domain $O$ such that
\begin{align*}
m\equiv \inf_{x\in O}V(x)<\min_{x\in \partial O}V(x).
\end{align*}
For the critical nonlinearity $f$, the similar results were obtained in  \cite{zjjzwm1,byeonzjj,zjjjmo}.

Now, we return our attention to  problem (\ref{question1}). In contrast to $s=1$, when $s \in(0, 1)$, $(-\DD)^s$ is a nonlocal operator, and some difficulties arise. Even in the subcritical case, there are only few references on the existence and concentration phenomena for (\ref{question1}).
In \cite{Davila}, J. D\'{a}vila, M. del Pino and J. Wei investigated (\ref{question1}) with $f(u)=u^p(1<p<2_s^*-1,2_s^*=2N/(N-2s))$. By applying the Lyapunov-Schmidt reduction method, they proved the existence of positive solutions which exhibit multiple spikes near given topologically nontrivial critical points of $V$ or cluster near a given local maximum point of $V$. More recently, C. O. Alves and O. H. Miyagaki \cite{Alvesfrac} considered (\ref{question1}) with a general nonlinearity $f$. By the  penalization method due to M. del Pino and P. Felmer \cite{M.delPino}, the authors constructed a spike solution around the local minimum point of $V$. In particular, in \cite{Alvesfrac}, the nonlinearity $f$ is subcritical and  satisfies the Ambrosetti-Rabinowitz condition and monotonicity condition. Inspired by \cite{Byeon1}, J. Seok \cite{seok} considered \re{question1} just under the Berestycki-Lions type conditions. By using the extension approach and local deformation argument,  J. Seok obtained positive solutions exhibiting multiple spikes near any given local minimum components of the potential $V$. More about the singularly perturbed fractional Schr\"{o}dinger problems, we refer to \cite{chenzheng,FMM,shangxzhangj} and the reference therein.

In the works above, only the subcritical case was considered. To study the semiclassical states of (\ref{question1}), the limit problem
\begin{align}
\label{limitequation2}
(-\Delta)^s u+mu=f(u),  u\in H^s(\mathbb{R}^N)
\end{align}
plays a crucial role. In many papers, the authors studied the existence of (ground state) solutions for fractional Schr\"{o}dinger equations when the nonlinearity $f$ satisfies the subcritical growth or critical growth (cf. \cite{bisci,tengkaiming,zjjjm}).
 For the critical case, the lack of compactness in the embedding $H^s(\mathbb{R}^N)\hookrightarrow L^{2_s^*}(\mathbb{R}^N)$ makes problem \re{question1} tough. A first breakthrough was given by X. He and W. Zou \cite{xiaominghe}. Precisely, the authors obtained the existence and concentration results for the problem $\varepsilon^{2s}(-\Delta)^su+V(x)u=g(u)+|u|^{2_s^*-2}u$. Here, we should point out that in \cite{xiaominghe}, $g$ satisfies the Ambrosetti-Rabinowitz condition and monotonicity condition.
A natural open question is that {\it in the critical case, whether a similar result as in \cite{xiaominghe} holds for a more general nonlinearity $f$, particularly without the Ambrosetti-Rabinowitz condition and monotonicity condition}. In this paper, we give an affirmative answer to this question.

\

\noindent{\bf 1.3 Main hypothesis.} In the present paper, we assume that $V$ satisfies $(V1)$-$(V2)$ and the nonlinearity $f$ satisfies

$(F_1)$  $f\in C^1(\mathbb{R}^+,\mathbb{R})$ and $\lim_{t\rightarrow 0}{f(t)/t}=0$,

$(F_2)$ $\lim_{t\rightarrow \infty}{f(t)/t^{2_s^*-1}}=1$,

$(F_3)$ there exist $\tilde{C}>0$ and $p<2_s^*$ such that $f(t)\geq t^{2_s^*-1}+\tilde{C}t^{p-1}$ for $t\geq 0$.

\

\noindent{\bf 1.4 Main result.}
\noindent Let
\begin{align*}
\mathcal{M}\equiv \{x\in O:V(x)=m\}.
\end{align*}

\noindent {\bf Theorem 1.1.} {\it Let $N>2s$ and $s\in (0,1)$. Suppose that $V\in C^1(\mathbb{R}^N,\mathbb{R})$ satisfies $(V_1)$-$(V_2)$
and $f$ satisfies $(F_1)$-$(F_3)$. Then for small $\varepsilon >0$, $(\ref{question1})$ admits a positive  solution $v_\varepsilon$ if $\max\{2_s^*-2,2\}<p<2_s^*$. Moreover, there exists a maximum point $y_\varepsilon$ of $v_\varepsilon$ such that $\lim_{\varepsilon\rightarrow 0}{\rm dist}(y_\varepsilon, \mathcal{M})=0$ and for any such $y_\varepsilon$, $w_\varepsilon(x)\equiv v_\varepsilon(\varepsilon x+y_\varepsilon)$ converges (up to a subsequence) uniformly to a least energy solution of $(\ref{limitequation2})$.
}

\

\noindent{\bf Remark 1.2.}   The condition $f\in C^1$ in  $(F_1)$ is to guarantee that a solution $u$ of (\ref{limitequation2}) satisfies the fractional Pohoz\v{a}ev identity. Since we are concerned with the positive solution, we assume $f(t)\equiv0$ for $t\leq 0$.

\

\noindent{\bf 1.5 Main difficulties and idea.}

The main difficulties are three-fold. Firstly, without the Ambrosetti-Rabinowitz condition, the boundedness of the (PS)-sequence is difficult to obtain. To overcome this difficulty, we seek the solutions in some neighborhood of
the set of ground state solutions to the limit problem.

Secondly, with the presence of the critical exponent $2_s^\ast$, the compactness of (PS)-sequence does not hold in general. To recover the compactness, we apply a truncation argument. Precisely, by the Moser iteration argument, we get a priori $L^\iy$-estimate of ground state solutions to the limit problem. Then we reduce the original problem to a subcritical problem and show the existence of spike solutions to the truncated problem. By the elliptic estimate, we show that the solution obtained is indeed a solution of the original problem.

Thirdly, in the truncation procedure, the uniform $L^\iy$-estimate of ground states to the limit problem plays a crucial role. However, the method introduced in \cite{B.Barriosa} can not be used directly. In this present paper, we prove that up to a translation, the set of ground state solutions is compact. By virtue of the compactness, we show that the ground state solutions are uniformly bounded in $L^\infty(\RN)$.

The paper is organized as follows. In Section 2, we introduce the variational setting and present some preliminary results. Section 3 is devoted to the study of ground state solutions of the limit problem (\ref{limitequation2}). The property of ground state solutions of (\ref{limitequation2}) such as  uniform boundedness is given by using the Moser iteration technique. Section 4 is devoted to the proof of Theorem 1.1.

\s {Preliminaries}
\renewcommand{\theequation}{2.\arabic{equation}}
By the scale change $x\rightarrow x/\varepsilon$ and denote $V_\varepsilon (x)=V(\varepsilon x)$,
then (\ref{question1}) is equivalent to
\begin{align}
\label{equivalentproblem}
(-\Delta)^s u+V_\varepsilon (x)u=f(u)\ \ \  \mbox{in} \ \ \mathbb{R}^N.
\end{align}
Thus, to study (\ref{question1}), it suffices to consider (\ref{equivalentproblem}). In the following, we present a quick survey of some preliminaries and properties about the fractional Sobolev spaces.

\

\noindent{\bf 2.1. Fractional  Sobolev spaces}

\noindent The fractional Laplacian operator $(-\Delta)^s$ with $s\in(0,1)$ of a function $u:\mathbb{R}^N\rightarrow \mathbb{R}$ is defined by
\begin{align}
\label{fourier}
\mathcal{F}((-\Delta)^s u)(\xi)=|\xi|^{2s}\mathcal{F}(u)(\xi),\ \ \xi\in\mathbb{R}^N,
\end{align}
where $\mathcal{F}$ is the Fourier transform. Consider the fractional Sobolev space
\begin{align*}
H^s(\mathbb{R}^N)=\{u\in L^2(\mathbb{R}^N):\int_{\mathbb{R}^N}(|\xi|^{2s}|\hat{u}|^2+\hat{u}^2)d\xi<\infty\},
\end{align*}
where $\hat{u}\doteq \mathcal{F}(u)$. The norm is defined by $\|u\|_{H^s(\mathbb{R}^N)}=(\int_{\mathbb{R}^N}(|\xi|^{2s}|\hat{u}|^2+\hat{u}^2)d\xi)^{\frac{1}2}$.
By Plancherel's theorem, we have $\|u\|_{L^2(\mathbb{R}^N)}=\|\hat{u}\|_{L^2(\mathbb{R}^N)}$ and for any $u\in H^s(\mathbb{R}^N)$,
\begin{align*}
\int_{\mathbb{R}^N}|(-\Delta)^{\frac{s}2}u(x)|^2dx=\int_{\mathbb{R}^N}|\widehat{(-\Delta)^{\frac{s}2}u(\xi)}|^2d\xi=\int_{\mathbb{R}^N}|\xi|^{2s}|\hat{u}|^2d\xi.
\end{align*}
It follows that $\|u\|_{H^s(\mathbb{R}^N)}=(\int_{\mathbb{R}^N}(|(-\Delta)^{\frac{s}2}u(x)|^2+u^2)dx)^{\frac{1}2}$, $u\in H^s(\mathbb{R}^N)$.

The space
$D^s(\mathbb{R}^N)$ is defined as the completion of $C_0^\infty(\mathbb{R}^N)$ under the norms
\begin{align*}
\|u\|_{D^s(\mathbb{R}^N)}^2=\int_{\mathbb{R}^N}|\xi|^{2s}|\hat{u}|^2d\xi=\int_{\mathbb{R}^N}|(-\Delta)^{\frac{s}2}u(x)|^2dx.
\end{align*}

\noindent Since we investigate problem (\ref{equivalentproblem}), we need the fractional Sobolev space  $H_{V_\varepsilon}^s(\mathbb{R}^N)$ which is a Hilbert space of $D^s(\mathbb{R}^N)$ with the norm
\begin{align*}
\|u\|_{H_{V_\varepsilon}^s(\mathbb{R}^N)}:=\left(\int_{\mathbb{R}^N}(|(-\Delta)^{\frac{s}2}u(x)|^2+V_\varepsilon(x) u^2)dx\right)^{\frac{1}2}<\infty.
\end{align*}

For the reader's convenience, we recall the embedding results for the fractional Sobolev spaces.

\

\noindent{\bf Lemma 2.1.}(see{\cite{PLLions}}) {\it Let $H_r^s(\mathbb{R}^N)=\{u\in H^s(\mathbb{R}^N):u(x)=u(|x|)\}$. For any $s\in(0,1), N>2s$, $H^s(\mathbb{R}^N)$ is continuously embedded into $L^q(\mathbb{R}^N)$ for $q\in [2,2_s^*]$ and compactly embedded into $L_{loc}^q(\mathbb{R}^N)$ for $q\in [1,2_s^*)$. Moreover, $H_r^s(\mathbb{R}^N)$ is compactly embedded into $L^q(\mathbb{R}^N)$ for $q\in (2,2_s^*)$.
}

\

\noindent{\bf Lemma 2.2.}(see{\cite{A.Cotsiolis,EGE}}) {\it For any $s\in(0,1)$, $D^s(\mathbb{R}^N)$ is continuously embedded into $L^{2_s^*}(\mathbb{R}^N)$, i.e., there exists $S_s>0$ such that
\begin{align*}
\|u\|_{L^{2_s^*}(\mathbb{R}^N)}\leq S_s \|u\|_{D^s(\mathbb{R}^N)}.
\end{align*}
}
\

\noindent{\bf 2.2. The variational setting}

\noindent Associated to (\ref{equivalentproblem}), the energy functional $I: H_{V_\varepsilon(x)}^s(\mathbb{R}^N)\rightarrow \mathbb{R}$ is defined by
\begin{align*}
I(u)=\frac{1}2\int_{\mathbb{R}^N}(|(-\Delta)^{\frac{s}2}u(x)|^2+V_\varepsilon(x)u^2)dx-\int_{\mathbb{R}^N}F(u)dx,\ \ \forall u\in H_{V_\varepsilon}^s(\mathbb{R}^N),
\end{align*}
where $F(t)=\int_0^t f(t)dt$. The conditions $(F_1)$-$(F_3)$ imply that $I(u)\in C^1$.

\

\noindent{\bf Definition 2.3.} We say that $u\in H_{V_\varepsilon}^s(\mathbb{R}^N)$ is a weak solution of (\ref{equivalentproblem}), if
\begin{align*}
\int_{\mathbb{R}^N} ((-\Delta)^{\frac{s}2}u (-\Delta)^{\frac{s}2}\phi+V_\varepsilon(x)u\phi) dx=\int_{\mathbb{R}^N} f(u)\phi dx, \ \forall \phi\in H_{V_\varepsilon}^s(\mathbb{R}^N).
\end{align*}

\noindent{\bf Proposition 2.4.}(See\cite{zjjjm}) {\it Under the assumptions of Theorem 1.1,
\begin{itemize}
\item [$(i)$] the limit problem $(\ref{limitequation2})$ admits a positive ground state solution, which is radially symmetric,

\item [$(ii)$]  let $S_m$ be the set of positive radial ground state solutions of $(\ref{limitequation2})$ whose maximum point is $0$, then  $S_m$ is compact in $H_r^s(\mathbb{R}^N)$,

\item [$(iii)$] for any solution $u\in S_m$, $u$ satisfies the Pohoz\v{a}ev identity
\begin{align}
\label{poho identity}
\frac{N-2s}2\int_{\mathbb{R}^N}|(-\Delta)^{\frac{s}2}u|^2=N\int_{\mathbb{R}^N}G(u),
\end{align}
where $G(u)=\int_{\mathbb{R}^N}F(u)-\frac{m}2u^2$,

\item [$(iv)$] for any solution $u\in S_m$, $u$ is also a mountain pass solution.
\end{itemize}
}
\

\

\noindent{\bf 2.3. Extended problems}

The fractional Laplacian operator is defined in the whole space through the Fourier transform(see(\ref{fourier})). In \cite{Caffarelli}, L. Caffarelli and L. Silvestre introduced that the fractional Laplacian operator also can be realized in a local way by using one more variable and the so-called $s$-harmonic extension. They developed a local interpretation of the fractional Laplacian operator given in $\mathbb{R}^N$ by considering a Dirichlet to Neumann type operator in the domain $\mathbb{R}_+^{N+1}:=\mathbb{R}^N\times(0,+\infty)$. More precisely, for $u\in D^s(\mathbb{R}^N)$, the solution $U\in D^1(t^{1-2s},\mathbb{R}_+^{N+1})$ of
\begin{align*}
\left\{
\begin{array}{ll}
-\mbox{div}(t^{1-2s}\nabla U)=0 \ \ \mbox{in} \ \ \ \mathbb{R}_+^{N+1},\\
U(x,0)=u \ \ \ \ \ \ \ \ \ \ \ \ \mbox{in} \ \ \mathbb{R}^N
\end{array}
\right.
\end{align*}
is called the $s$-harmonic extension of $u$ and satisfies
\begin{align*}
-\lim_{t\rightarrow 0}t^{1-2s}\partial_t U(x,t)=N_s(-\Delta)^su(x)\ \ \mbox{in}\ \ \ \mathbb{R}^N,
\end{align*}
where $N_s=2^{1-2s}\Gamma(1-s)/\Gamma(s)$. We denotes $U\doteq E_s(u)$ and $u=tr(U)=U(x,0)$. From \cite{Caffarelli}, the  $s$-harmonic extension of $u$ is defined by
\begin{align}
\label{harmonic expression}
U(x,t)=\int_{\mathbb{R}^N}P_s(x-\xi,t)u(\xi)d\xi,
\end{align}
where $P_s(x,t)=C_{N,s}\frac{t^{2s}}{\left(|x|^2+|t|^2\right)^{\frac{N+2s}2}}$ with the constant $C_{N,s}$ satisfying $\int_{\mathbb{R}^N}P_s(x,1)dx=1$.

The space $D^1(t^{1-2s},\mathbb{R}_+^{N+1})$ denotes the completion of $C_0^\infty(\overline{\mathbb{R}_+^{N+1}})$ with the norm
\begin{align*}
\|U\|_{D^1(t^{1-2s},\mathbb{R}_+^{N+1})}^2=\int_{\mathbb{R}_+^{N+1}}t^{1-2s}|\nabla U(x,t)|^2dxdt,
\end{align*}
which satisfies(see\cite{c.brandle})
\begin{align}
\label{extennorm}
\|U\|_{D^1(t^{1-2s},\mathbb{R}_+^{N+1})}=\sqrt{N_s}\|u\|_{D^s(\mathbb{R}^N)}.
\end{align}
It is known that(see\cite{jinlianling}) for any $U\in D^1(t^{1-2s},\mathbb{R}_+^{N+1})$, the trace $U(x,0)$ belongs to $D^s(\mathbb{R}^N)$ and the trace map is continuous as follows,
\begin{align}
\|U(x,0)\|_{D^s(\mathbb{R}^N)}\leq C\|U\|_{D^1(t^{1-2s},\mathbb{R}_+^{N+1})}.
\end{align}

By  the $s$-harmonic extension, we introduce extended problems of (\ref{limitequation2}) and (\ref{equivalentproblem}) respectively,
\begin{align}
\label{extendedpro}
\left\{
\begin{array}{ll}
-\mbox{div}(t^{1-2s}\nabla U)=0 \ \ \mbox{in} \ \ \ \mathbb{R}_+^{N+1},\\
-\frac{1}{N_s}\lim_{t\rightarrow 0}t^{1-2s}\partial_t U(x,t)=-mU(x,0)+f(U(x,0)) \ \mbox{in} \ \mathbb{R}^N
\end{array}
\right.
\end{align}
and
\begin{align}
\label{equalextendedpro}
\left\{
\begin{array}{ll}
-\mbox{div}(t^{1-2s}\nabla U)=0 \ \ \mbox{in} \ \ \ \mathbb{R}_+^{N+1},\\
-\frac{1}{N_s}\lim_{t\rightarrow 0}t^{1-2s}\partial_t U(x,t)=-V_\varepsilon (x)U(x,0)+f(U(x,0)) \ \mbox{in} \ \mathbb{R}^N.
\end{array}
\right.
\end{align}

Define the function spaces $H_0$ and $H_{V_\varepsilon}$ by the sets of $U\in D^1(t^{1-2s},\mathbb{R}_+^{N+1})$ satisfying
\begin{align*}
&\|U\|_0^2=\int_{\mathbb{R}_+^{N+1}}t^{1-2s}|\nabla U(x,t)|^2dxdt+\int_{\mathbb{R}^N}U^2(x,0)dx<\infty,\\
&\|U\|_{V_\varepsilon}^2=\int_{\mathbb{R}_+^{N+1}}t^{1-2s}|\nabla U(x,t)|^2dxdt+\int_{\mathbb{R}^N}V_\varepsilon(x)U^2(x,0)dx<\infty.
\end{align*}
We call $U\in H_0$ a weak solution of (\ref{extendedpro}) if $U$ satisfies
\begin{align*}
\int_{\mathbb{R}_+^{N+1}}t^{1-2s}\nabla U\nabla Vdxdt+N_s\int_{\mathbb{R}^N}[mU(x,0)-f(U(x,0))]V(x,0)dx=0
\end{align*}
for any $V\in H_0$. It is well known that if $U\in H_0$ is a weak solution of (\ref{extendedpro}), then $U(x,0)\in H^s(\mathbb{R}^N)$ is a weak solution of (\ref{limitequation2}). Similarly we can define the weak solution of (\ref{equalextendedpro}) and get the similar relation between solutions of (\ref{equivalentproblem}) and (\ref{equalextendedpro}).

Proposition 2.4 gives  the main results about the existence and compactness of ground state solutions of (\ref{limitequation2}). In the following, we investigate the existence and compactness of ground state solutions of (\ref{extendedpro}).

\

\noindent{\bf Proposition 2.5.} {\it Under the assumptions of Theorem 1.1, there exists a positive ground solution to $(\ref{extendedpro})$. Denote by $\tilde{S}_m$ the set of ground state solutions of $(\ref{extendedpro})$ such that, for any $U\in \tilde{S}_m$, $tr(U)$ is radial and attains its maximum at $0\in \mathbb{R}^N$, then $\tilde{S}_m$ is compact.
}
\

\bp
From Proposition 2.4, we know that under the assumptions of Theorem 1.1, there exists a least energy solution of (\ref{limitequation2}). Let $S_m$ be the set of positive radial and decreasing ground state solutions, then  $S_m$ is compact in $H_r^s(\mathbb{R}^N)$.
Define the energy functional associated to (\ref{limitequation2}) and (\ref{extendedpro}) respectively by
\begin{align*}
&I(u)=\frac{1}2\int_{\mathbb{R}^N}|(-\Delta)^{\frac{s}2}u|^2dx+\int_{\mathbb{R}^N}[\frac{1}2 mu^2-F(u)]dx,\\
&J(U)=\frac{1}2\int_{\mathbb{R}_+^{N+1}}t^{1-2s}|\nabla U|^2dxdt+N_s\int_{\mathbb{R}^N}[\frac{1}2 mU(x,0)^2-F(U(x,0))]dx.
\end{align*}
From the argument about extended problems, by a simple calculation, we can obtain that $J(U)=N_sI(u)$. So we can see that the solutions of (\ref{limitequation2}) and (\ref{extendedpro}) are one to one.

Let $\{u_n\}\in S_m$, then up to a subsequence, $u_n\rightarrow u_0$ strongly in $S_m$ since $S_m$ is compact. Denote by $U_n=E_s(u_n)$ and $U_0=E_s(u_0)$ the solutions of (\ref{extendedpro}) corresponding to $u_n$ and  $u_0$ respectively. Then $U_n\in \tilde {S}_m$ and $U_0\in \tilde {S}_m$. To prove $\tilde {S}_m$ is compact, it suffices to prove that $\|U_n-U_0\|_0\rightarrow 0$. Since $U_n$ and $U_0$ are (weak) solutions of (\ref{extendedpro}), for any $\phi\in H_0$, we have
\begin{align}
\label{weaksolution1}
\int_{\mathbb{R}_+^{N+1}}t^{1-2s}(\nabla U_n\nabla \phi)dxdt
=N_s\int_{\mathbb{R}^N}[f(U_n(x,0)-mU_n(x,0)]\phi(x,0)dx
\end{align}
and
\begin{align}
\label{weaksolution2}
\int_{\mathbb{R}_+^{N+1}}t^{1-2s}(\nabla U_0\nabla \phi)dxdt
=N_s\int_{\mathbb{R}^N}[f(U_0(x,0)-mU_0(x,0)]\phi(x,0)dx.
\end{align}
From (\ref{weaksolution1})-(\ref{weaksolution2}), we get
\begin{align*}
&\int_{\mathbb{R}_+^{N+1}}t^{1-2s}(\nabla U_n-\nabla U_0)\nabla \phi dxdt\\
&=N_s\int_{\mathbb{R}^N}[f(U_n(x,0))-f(U_0(x,0))-m(U_n(x,0)-U_0(x,0))]\phi(x,0)dx.
\end{align*}

\noindent Since $u_n\rightarrow u_0$ in $H^s(\mathbb{R}^N)$, we have
\begin{align*}
\lim_{n\rightarrow \infty}\int_{\mathbb{R}_+^{N+1}}t^{1-2s}(\nabla U_n-\nabla U_0)\nabla \phi=0.
\end{align*}
Let $\phi=U_n-U_0$, then
\begin{align*}
\lim_{n\rightarrow \infty}\int_{\mathbb{R}_+^{N+1}}t^{1-2s}|\nabla U_n-\nabla U_0|^2=0,
\end{align*}
which implies $\lim_{n\rightarrow \infty}\|U_n-U_0\|_{D^1(t^{1-2s},\mathbb{R}_+^{N+1})}^2=0$.
Let $n\rightarrow \infty$, we can obtain
\begin{align*}
\|U_n-U_0\|_{0}^2=\|U_n-U_0\|_{D^1(t^{1-2s},\mathbb{R}_+^{N+1})}^2+\|U_n(x,0)-U_0(x,0)\|_{L^2(\mathbb{R}^N)}^2\rightarrow 0.
\end{align*}
The proof is completed.
\ep
\

For studying the existence and concentration phenomena by harmonic extension, we need some elliptic estimates for extended problems.

\

\noindent{\bf 2.4. Elliptic estimates}

\noindent Let $\Omega_r:=B_r^N(0)\times (0,r)$. Consider the following nonlinear Neumann boundary value problem
\begin{align}
\label{ellipticestimateequation}
\left\{
\begin{array}{ll}
-\mbox{div}(t^{1-2s}\nabla U)=0 \ \ \mbox{in} \ \ \Omega_1,\\
-\lim_{t\rightarrow 0}t^{1-2s}\partial_t U(x,t)=a(x)U(x,0)+g(x) \ \mbox{in} \ B_1^N(0).
\end{array}
\right.
\end{align}
Let $H^1(t^{1-2s},\Omega_r)$ be the weighted Sobolev space with the norm
\begin{align*}
\|U(x,t)\|_{H^1(t^{1-2s},\Omega_r)}^2&=\int_{\Omega_r}|t|^{1-2s}(U^2+|\nabla U|^2)\\
&=\|U(x,t)\|_{L^2(t^{1-2s},\Omega_r)}^2+\|U(x,t)\|_{D^1(t^{1-2s},\Omega_r)}^2.
\end{align*}

\

\noindent{\bf Proposition 2.6.}(De Giorgi-Nash-Moser type estimate.)(see\cite{jinlianling,seok}) {\it Suppose $a,g\in L^p(B_1^N(0))$ for some $p>\frac{N}{2s}$.

(i)Let $U\in H^1(t^{1-2s},\Omega_1)$ be a weak solution of $(\ref{ellipticestimateequation})$. Then $U\in L^\infty(\Omega_{1/2})$ and there exist a constant $C>0$ depending only on $N,s,p$ and $\|a\|_{L^p(B_1^N(0))}$ such that
\begin{align}
\label{ellipticestimate1}
\sup_{\Omega_{1/2}}U\leq C\left(\|U\|_{L^2(t^{1-2s},\Omega_1)}+\|g\|_{L^p(B_1^N(0))}\right);
\end{align}

(ii)Let $U\in H^1(t^{1-2s},\Omega_1)$ be a weak solution of $(\ref{ellipticestimateequation})$. Then there is a $\alpha\in (0,1)$ depending only on $N,s,p$ such that $U\in C^\alpha(\overline{\Omega_{1/2}})$
and there exist a constant $C>0$ depending only on $N,s,p$ and $\|a^+\|_{L^p(B_1^N(0))}$ such that
\begin{align}
\label{ellipticestimate2}
\|U\|_{C^\alpha(\overline{\Omega_{1/2}})}\leq C\left(\|U^+\|_{L^\infty(\Omega_1)}+\|g\|_{L^p(B_1^N(0))}\right).
\end{align}
}

\s{A priori estimate of ground state solutions for limit problems (\ref{limitequation2}) and (\ref{extendedpro})}
\renewcommand{\theequation}{3.\arabic{equation}}
In this section, we study a priori $L^\infty$-estimate of ground state solutions to the limit problems (\ref{limitequation2}) and (\ref{extendedpro}). However, the method introduced in \cite{B.Barriosa} is only used to get the $L^\infty$-estimate for any fixed solution. With the help of the compactness of $S_m$, we modify the argument in \cite{B.Barriosa} to get the uniform $L^\infty$-boundedness of ground state solutions to (\ref{limitequation2}) and (\ref{extendedpro}).

\

\noindent{\bf Proposition 3.1.} {\it Under the assumptions of Theorem 1.1, for any $u\in S_m$, $u\in L^\infty(\mathbb{R}^N)$. Moreover, $\sup\{\|u\|_\infty:u\in S_m\}<\infty.$}

\bp
It suffices to prove that for any $\{u_n\}\subset S_m$ with $u_n\rg u_0$ strongly in $H^s(\RN)$, $u_n\in L^\infty(\mathbb{R}^N)$ and $\sup_n\|u_n\|_\infty<\infty.$ Let $\gamma\geq 1$, $T>0$, we define
\begin{align*}
\varphi(t)=\varphi_{\gamma,T}(t)=\left\{
\begin{array}{ll}
0,\ \ \ \ \ \ \ \ \ \ \ \ \ \ \ \ \ \ \ \ \ \ \ \ \ \ t\leq 0,\\
t^{\gamma},\ \ \ \ \ \ \ \ \ \ \ \ \ \ \ \ \ \ \ \ \ \ \ \ \ 0<t<T,\\
\gamma T^{\gamma-1}(t-T)+T^{\gamma},\ \ \ t\geq T.
\end{array}
\right.
\end{align*}
It is easy to verify that, for any $t\in\mathbb{R}$, $\varphi'(t)\leq \gamma T^{\gamma-1}$ and $t\varphi'(t)\leq \gamma\varphi(t)$. Since $\varphi(t)$ is convex, we have $(-\Delta)^s\varphi(u_n)\leq \varphi'(u_n)(-\Delta)^s u_n$ and consequently $\|\varphi(u_n)\|_{D^s(\mathbb{R}^N)}\leq \gamma T^{\gamma-1}\|u_n\|_{D^s(\mathbb{R}^N)}$. On the other hand, from Lemma 2.2, we obtain
\begin{align}
\label{guji1}
S_s\|\varphi(u_n)\|_{D^s(\mathbb{R}^N)}\geq \|\varphi(u_n)\|_{L^{2_s^*}(\mathbb{R}^N)}.
\end{align}
It follows from $(F_1)$ and $(F_2)$ that there exists a constant $C$ such that $f(t)\leq \frac{mt}2+Ct^{2_s^*-1}$ for $t>0$. Noting that  $u_n\geq 0$, we have $(-\DD)^su_n\leq Cu_n^{2_s^*-1}$ in $\RN$ and then
\begin{align}
\label{guji2}
\|\varphi(u_n)\|_{D^s(\mathbb{R}^N)}^2=\int_{\mathbb{R}^N} \varphi(u_n)\varphi'(u_n)(-\Delta)^su_n\leq C\int_{\mathbb{R}^N} \varphi(u_n)\varphi'(u_n)u_n^{2_s^*-1}.
\end{align}
Together with (\ref{guji1}) and $u_n\varphi'(u_n)\leq \gamma\varphi(u_n)$, we deduce that
\begin{align}
\label{guji3}
\|\varphi(u_n)\|_{L^{2_s^*}(\mathbb{R}^N)}^2\leq CS_s^2\int_{\mathbb{R}^N}\varphi'(u_n)u_n \varphi(u_n)u_n^{2_s^*-2}\leq C_\gamma \int_{\mathbb{R}^N} \varphi^2(u_n)u_n^{2_s^*-2},
\end{align}
where $C_\gamma=CS_s^2\gamma$. By H\"{o}lder's inequality, we infer that
\begin{align}
\label{guji4}
\int_{\mathbb{R}^N}\varphi^2(u_n)u_n^{2_s^*-2}=\int_{\{u_n\leq R\}}\varphi^2(u_n)u_n^{2_s^*-2}+\int_{\{u_n> R\}}\varphi^2(u_n)u_n^{2_s^*-2}\nonumber\\
\leq \int_{\{u_n\leq R\}}\varphi^2(u_n)R^{2_s^*-2}+\|\varphi(u_n)\|_{L^{2_s^*}(\mathbb{R}^N)}^2\left(\int_{\{u_n> R\}}u_n^{2_s^*}\right)^{\frac{2_s^*-2}{2_s^*}}.
\end{align}
Since $u_n\rg u_0$ strongly in $H^s(\RN)$, we take $R$ large enough(be fixed later) such that $\left(\int_{\{u_n> R\}}u_n^{2_s^*}\right)^{\frac{2_s^*-2}{2_s^*}}\leq \frac{1}{2C_\gamma}.$ Using the fact that $\varphi(u_n)\leq u_n^{\gamma}$, it follows from (\ref{guji3}) and (\ref{guji4}) that
\begin{align}
\label{guji5}
\|\varphi(u_n)\|_{L^{2_s^*}(\mathbb{R}^N)}^2\leq \tilde{C_\gamma}R^{2_s^*-2}\int_{\mathbb{R}^N}u_n^{2\gamma},
\end{align}
where $\tilde{C_\gamma}>0$ is a constant. If $u_n\in L^{2\g}(\RN)$, letting $T\rightarrow \infty$, we  get $u_n\in L^{2_s^*\gamma}(\mathbb{R}^N).$ By iteration, for any $p\geq 2$, $u_n\in L^p(\RN)$.
In (\ref{guji3}), letting $T\rightarrow \infty$ and using the fact $\varphi(u_n)\leq u_n^{\gamma}$, we infer that $\|u_n\|_{L^{2_s^*\gamma}(\mathbb{R}^N)}^{2\gamma}\leq C_\gamma\int_{\mathbb{R}^N}u_n^{2_s^*+2\gamma-2}$.

Let $\g_1=2_s^*/2$ and $2_s^*+2\g_{i+1}-2=2_s^*\g_i, i=1,2,\cdots$, then
\begin{align*}
\gamma_{i+1}-1=\left(\frac{2_s^*}2\right)^i(\gamma_1-1)
\end{align*}
and
\begin{align}
\label{jielun1}
\left(\int_{\mathbb{R}^N}u_n^{\gamma_{i+1}2_s^*}\right)^{\frac{1}{2_s^*(\gamma_{i+1}-1)}}\leq C_{\gamma_{i+1}}^{\frac{1}{2(\gamma_{i+1}-1)}}\left(\int_{\mathbb{R}^N}u_n^{\gamma_i2_s^*}\right)^{\frac{1}{2_s^*(\gamma_{i}-1)}},
\end{align}
where $C_{\gamma_{i+1}}=CS_s^2\g_{i+1}$. Denote $K_{i}=\left(\int_{\mathbb{R}^N}u_n^{\gamma_i2_s^*}\right)^{\frac{1}{2_s^*(\gamma_{i}-1)}}$, then for $\tau>0$,
\begin{align*}
K_{\tau+1}\leq \prod_{i=2}^{\tau}C_{\gamma_{i}}^{\frac{1}{2(\gamma_{i}-1)}}K_1.
\end{align*}
After a simple calculation, we conclude that there exists a constant $C_0>0$, independent of $\tau$ and $n$, such that $K_{\tau+1}\leq C_0K_1$ for any $n\in N$.
Hence,
\begin{align*}
\|u_n\|_{L^\infty(\mathbb{R}^N)}\leq C_0K_1<\infty, \ \ \forall n\in N.
\end{align*}
The proof is finished.
\ep

\

%
%
%

\s{Proof of Theorem 1.1}
\renewcommand{\theequation}{4.\arabic{equation}}
In this section, we use the truncation approach to prove Theorem 1.1. First, we construct spike solutions of the truncation problem in some neighborhood of ground sate solutions to the limit problem. Second, by the elliptic estimate we show that a solution of the truncation problem is indeed a solution of the original problem.

\subsection{Truncation problems}
By Proposition 3.1, there exists $\bar{k}>0$ such that
\begin{align*}
\sup _{u\in S_m}\|u\|_{\iy}<\bar{k}.
\end{align*}
For any $k>\max_{t\in[0,\bar{k}]}f(t)$, let
$
f_k(t)=\min\{f(t),k\}(t\in\mathbb{R}).
$
In the following, we consider the truncation problem
\begin{align}
\label{trunctionfk}
(-\Delta)^su+V_\varepsilon(x)u=f_k(u), \ \ \ u\in H_{V_\varepsilon}^s(\mathbb{R}^N).
\end{align}
Obviously, any solution $u_\e$ of \re{trunctionfk} is indeed a solution of the original problem \re{equivalentproblem} if $\|u_\e\|_\iy\le \bar{k}$.
Now, we consider the corresponding limit equation
\begin{align}
\label{trunctionequation}
(-\Delta)^su+mu=f_k(u), \ \ \ u\in H^s(\mathbb{R}^N),
\end{align}
whose energy functional is given by
\begin{align*}
I_m^k(u)=\frac{1}2\int_{\mathbb{R}^N}|(-\Delta)^{\frac{s}2}u|^2+mu^2-\int_{\mathbb{R}^N}F_k(u),
\end{align*}
and the extended problem is
\begin{align}
\label{trunctionextendedpro}
\left\{
\begin{array}{ll}
-\mbox{div}(t^{1-2s}\nabla U)=0 \ \ \mbox{in} \ \ \ \mathbb{R}_+^{N+1},\\
-\frac{1}{N_s}\lim_{t\rightarrow 0}t^{1-2s}\partial_t U(x,t)=-mU(x,0)+f_k(U(x,0)) \ \mbox{in} \ \mathbb{R}^N,
\end{array}
\right.
\end{align}
where $F_k(t)=\int_0^tf_k(t)dt$.

\

\noindent{\bf Lemma 4.1.} For any $k>\max_{t\in[0,\bar{k}]}$, (\ref{trunctionextendedpro}) admits a positive least energy solution.
\bp It suffices to prove $f_k$ satisfies the Berestycki-Lions type conditions.
It is obvious that $f_k(t)=o(t)$ as $t\rightarrow 0$ and $\limsup_{t\rightarrow \infty}{f_k(t)/t^p}<C$ for some $C>0$ and $p\in (1,2_s^*-1)$. Now, we show that there exists $T>0$ such that $mT^2<2F_k(T)$.
Indeed, taking any $u\in S_m$, by the Pohoz\v{a}ev identity
\begin{align}
\label{poho identity}
\frac{N-2s}2\int_{\mathbb{R}^N}|(-\Delta)^{\frac{s}2}u|^2=N\int_{\mathbb{R}^N}F(u)-\frac{m}2u^2,
\end{align}
which means that there exists $x_0\in \mathbb{R}^N$ such that $F(u(x_0))>\frac{m}2u^2(x_0)$. Obviously, $F_k(u(x))\equiv F(u(x))$ for all $x\in \mathbb{R}^N$. Let $T=u(x_0)>0$, we get that $F_k(T)>\frac{m}2T^2$.
\ep

As can be seen in J. Seok\cite{seok}, for every least energy solution $U(x,t)$ of (\ref{trunctionextendedpro}), the trace $u=U(x,0)$ is a  positive classical solution of (\ref{trunctionequation}). Moreover, $U(x,t)$ is a mountain pass solution. Denote by $\tilde{S}_m^k$ the set of least energy solutions $U$ of (\ref{trunctionextendedpro}) such that $U(x,0)$ attains its maximum at $0\in\mathbb{R}^N$, then $\tilde{S}_m^k$ is compact. Denote by $E_m^k$ the energy of $u\in S_m^k$. Noting that $f_k(t)\leq f(t)$ for any $t$, we get that $E_m^k\geq E_m$ because  every solution is a mountain pass solution. On the other hand, it follows from $\sup _{u\in S_m}\|u\|_{\iy}<\bar{k}$ and the definition of $f_k$ that $S_m\subset S_m^k$. So $E_m^k\leq E_m$. Thus,
\begin{align*}
E_m^k=E_m\ \ \ \ \mbox{for} \ \ \ k>\max_{t\in [0,\bar{k}]}f(t).
\end{align*}

\

\noindent{\bf Lemma 4.2.} {\it For $k>\max_{t\in [0,\bar{k}]}f(t)$, we have
$
S_m^k=S_m.
$

\bp Obviously, $S_m\subset S_m^k$. Now, we claim that $S_m^k\subset S_m$. Let
\begin{align*}
&T(u)=\int_{\mathbb{R}^N}|(-\Delta)^{\frac{s}{2}}u|^2,\,\, V(u)=\int_{\mathbb{R}^N}G(u),\,\, V_k(u)=\int_{\mathbb{R}^N}G_k(u),
\end{align*}
where $G(u)=F(u)-\frac{m}2u^2,\,\, G_k(u)=F_k(u)-\frac{m}2u^2.$ We consider the constraint minimization problems
\begin{align}
\label{constraintproblem1}
M:=\inf\{T(u):V(u)\equiv1,u\in H^s(\mathbb{R}^N)\}
\end{align}
and
\begin{align}
\label{constraintproblem2}
M_k:=\inf\{T(u):V_k(u)\equiv1,u\in H^s(\mathbb{R}^N)\}.
\end{align}
Let $\tilde{u}_k$ is a minimizer of (\ref{constraintproblem2}), then by the Pohoz\v{a}ev's identity (\ref{poho identity}), $u_k=\tilde{u}_k(\frac{x}\sigma)\in S_m^k$, where $\sigma=(\frac{N-2s}{2N}M_k)^{\frac{1}{2s}}$. Moreover, $u_k$ is a minimizer of $T(u)$ in
\begin{align*}
\left\{u\in H^s(\mathbb{R}^N):V_k(u)=\left(\frac{N-2s}{2N}M_k\right)^{\frac{N}{2s}}\right\}
\end{align*}
and
\begin{align*}
E_m^k=I_m^k(u_k)=\frac{1}2T(u_k)-V_k(u_k)=\frac{s}N\left(\frac{N-2s}{2N}\right)^{\frac{N-2s}{2s}}M_k^{\frac{N}{2s}}.
\end{align*}
Similarly, we obtain $E_m=\frac{s}N\left(\frac{N-2s}{2N}\right)^{\frac{N-2s}{2s}}M^{\frac{N}{2s}}$. Since $E_m=E_m^k$, we have $M=M_k$ for $k>\max_{t\in [0,\bar{k}]}f(t)$. Now, we claim that $\tilde{u}_k$ is also a minimizer of (\ref{constraintproblem1}). As a consequence, $u_k=\tilde{u}_k(\frac{x}\sigma)\in S_m$ which means
$S_m^k\subset S_m$. Since $\tilde{u}_k$ is a minimizer of (\ref{constraintproblem2}), we have
\begin{align}
\label{tukmk1}
T(\tilde{u}_k)=M_k=M,
\end{align}
and
\begin{align}
\label{tukmk2}
T(\tilde{u}_k)=M_kV_k(\tilde{u}_k)=M_k\left(V_k(\tilde{u}_k)\right)^{\frac{N-2s}N}\leq M\left(V(\tilde{u}_k)\right)^{\frac{N-2s}N}.
\end{align}
Now, it suffices to show that
$
V(\tilde{u}_k)=1.
$
Let $\tilde{w}_k=\tilde{u}_k(\lambda\cdot)$ such that $V(\tilde{w}_k)=\lambda^{-N}V(\tilde{u}_k)=1$, where $\lambda=\left(V(\tilde{u}_k)\right)^{\frac{1}N}$. Thus we get $T(\tilde{w}_k)=\lambda^{-N+2s}T(\tilde{u}_k)=\left(V(\tilde{u}_k)\right)^{\frac{-N+2s}N}T(\tilde{u}_k)\geq M$.
So we have
$
T(\tilde{u}_k)\geq M\left(V(\tilde{u}_k)\right)^{\frac{N-2s}N}.
$
Together with (\ref{tukmk2}), we obtain
$
T(\tilde{u}_k)=M\left(V(\tilde{u}_k)\right)^{\frac{N-2s}N}.
$
It follows from (\ref{tukmk1}) that $V(\tilde{u}_k)=1$. The proof is completed.
\ep

\noindent{\bf Corollary 4.3.} {\it For $k>\max_{t\in [0,\bar{k}]}f(t)$,
$
\tilde{S}_m^k=\tilde{S}_m.
$
\br
By Proposition 2.5, up to translations, $\tilde{S}_m^k$ is compact.
\er
\subsection{Completion of Proof of Theorem 1.1.}
\bp
{\bf Step 1.} The existence of solutions to the truncation problem (\ref{trunctionfk}).
By Corollary 4.3, for fixed $k>\max_{t\in [0,\bar{k}]}f(t)$, we have $\tilde{S}_m^k=\tilde{S}_m$. In the following, we construct a solution of (\ref{trunctionfk}) in some neighborhood of $\tilde{S}_m$. Precisely, define a set of approximation solutions by
\begin{align*}
N_\varepsilon(\rho)=\{\phi_\varepsilon(\cdot-\frac{x_\varepsilon}\varepsilon)U(\cdot-\frac{x_\varepsilon}\varepsilon,t)+w:x_\varepsilon\in M^\delta,U\in\tilde{S}_m,\\ w\in H_{V_\varepsilon},\delta>0,\|w\|_{V_\varepsilon}\leq \rho\},
\end{align*}
where $\phi:\mathbb{R}^N\rightarrow[0,1]$ is defined as a smooth cut-off function satisfying
\begin{align*}
\phi(x)=\left\{
\begin{array}{ll}
1,\ \ \ |x|\leq \delta,\\
0, \ \ \ |x|\geq 2\delta,
\end{array}
\right.
\end{align*}
and $M^\delta=\{x\in \mathbb{R}^N:dist(x,A)\leq \delta\}$. By Lemma 4.1, $f_k$ satisfies all the assumptions in \cite[Theorem 1.1]{seok}, which implies that (\ref{trunctionfk}) admits a solution $u_\e$, where $u_\varepsilon=U_\varepsilon(\cdot,0)$ and $U_\varepsilon\in N_\varepsilon(\rho)$ is a solution of (\ref{equalextendedpro}) where $f(U(x,0))$ is replaced by $f_k(U(x,0))$. Moreover, there exists a maximum point $x_\varepsilon$ of $u_\varepsilon$ such that $\lim_{\varepsilon\rightarrow 0}dist(\varepsilon x_\varepsilon,\mathcal{M})=0$. As a consequence of \cite[Proposition 3.1]{seok}, $\|u_\varepsilon(\cdot+x_\varepsilon)\rightarrow u(\cdot +z_0)\|_{H_{V_\varepsilon}^s(\mathbb{R}^N)}\rg0$ as $\varepsilon\rightarrow 0$, where $u\in S_m^k=S_m$  and  $z_0\in \mathbb{R}^N$.

{\bf Step 2.} We show that $u_\e$ is indeed a solution of the original problem (\ref{equivalentproblem}). By Step 1, $u_\varepsilon(\cdot+x_\varepsilon)\rightarrow u(\cdot +z_0)$ strongly in $H^s(\mathbb{R}^N)$ as $\varepsilon\rightarrow 0$. Then, similar as in Proposition 3.1, we know that there exists $\varepsilon_0>0$, $\sup_{\e\leq\e_0}\|u_\varepsilon\|_\iy<\iy$. If follows from (\ref{harmonic expression}) and the fact $\int_{\mathbb{R}^N}P_s(x,1)dx=1$, that for any $x\in\RN, t\in\R^+$,
\begin{align*}
|U_\e(x,t)|\leq \|u_\e\|_{L^\infty(\mathbb{R}^N)}\int_{\mathbb{R}^N}P_s(x,t)dx=\|u_\e\|_{L^\infty(\mathbb{R}^N)}.
\end{align*}
So $\sup_{\e\leq\e_0}\|U_\varepsilon\|_{L^\infty(\mathbb{R}^N\times\R^+)}<\iy$. Let $w_\varepsilon(\cdot)=u_\varepsilon(\cdot+x_\varepsilon)$, then $w_\varepsilon(0)=\|u_\e\|_\iy$. To conclude the proof, it suffices to show that for $\e>0$ small enough, $w_\e(0)<\bar{k}$. Indeed, $w_\varepsilon=U_\e(x+x_\e,0)$ and $U_\e(x+x_\e,t)$ satisfies
\begin{align*}
\left\{
\begin{array}{ll}
-\mbox{div}(t^{1-2s}\nabla U_\e(x+x_\e,t))=0 \ \ \mbox{in} \ \ \ \mathbb{R}_+^{N+1},\\
-\frac{1}{N_s}\lim_{t\rightarrow 0}t^{1-2s}\partial_t U_\e(x+x_\e,t)=-V_\varepsilon (x+x_\e)U_\e(x+x_\e,0)+f_k(U_\e(x+x_\e,0)) \ \mbox{in} \ \mathbb{R}^N.
\end{array}
\right.
\end{align*}
Then by virtue of Proposition 2.6, for some $\alpha\in (0,1)$ and $C>0$ (independent of $\e$), such that
$$
\|U(\cdot+x_\e,\cdot)\|_{C^\alpha(\overline{\Omega_{1/2}})}\leq C\left(\|U_\e\|_{L^\infty(\RN\times\R^+)}+\|f_k(U_\e(x+x_\e,0))\|_{L^p(B_1^N(0))}\right),
$$
where $\overline{\Omega_{1/2}}=\overline{B_{1/2}^N(0)}\times [0,1/2]$. It follows from $\sup_{\e\leq\e_0}\|U_\varepsilon\|_{L^\infty(\mathbb{R}^N\times\R^+)}<\iy$ that
$$
\sup_{\e\leq\e_0}\|w_\e\|_{C^\al(B_{1/2}^N(0))}<\iy.
$$
It implies that $\{w_\varepsilon\}$ is uniformly bounded and equicontinuous concerning $\varepsilon$ in $B_1^N(0)$. By  Arzel\`{a}-Ascoli theorem,  $w_\varepsilon(\cdot)\rightarrow u(\cdot +z_0)$ uniformly in $B_1^N(0)$ and then $\|u_\varepsilon\|_\iy=w_\varepsilon(0)<\bar{k}$ uniformly holds for sufficiently small $\varepsilon>0$. Therefore, for $\e>0$ small, $f_k(u_\varepsilon(x))\equiv f(u_\varepsilon(x))$ for $x\in\RN$, which means that $u_\varepsilon(x)$ is a solution of the original problem (\ref{equivalentproblem}). Let $v_\varepsilon(\cdot)=u_\varepsilon(\cdot/\varepsilon)$ and $y_\varepsilon=\varepsilon x_\varepsilon$, then $v_\varepsilon$ is a solution of (\ref{question1}), whose maximum point is $y_\varepsilon$ satisfying $\lim_{\varepsilon\rightarrow 0}\mbox{dist}(y_\varepsilon,\mathcal{M})=0$. The proof is completed.
\ep

\noindent{\bf Acknowledgements.}\,\,
{\rm The first and second authors were supported  by the National Natural Science Foundation of China (11271364). The third author was supported by the Science Foundation of Chongqing Jiaotong University (15JDKJC-B033).}

\end{document}